\newif\ifJOURNAL
\JOURNALfalse
\newif\ifWP
\WPfalse
\newif\ifBASIC
\BASICfalse

\newif\ifFULL
\FULLfalse
\newif\ifLATIN
\LATINfalse

\WPtrue

\LATINtrue	

\newif\ifnotJOURNAL		
\notJOURNALtrue
\ifJOURNAL\notJOURNALfalse\fi

\newif\ifnotWP		
\notWPtrue
\ifWP\notWPfalse\fi

\newif\ifnotFULL	
\notFULLtrue
\ifFULL\notFULLfalse\fi

\newif\ifnotLATIN	
\notLATINtrue
\ifLATIN\notLATINfalse\fi

\ifnotLATIN

\fi
\ifLATIN

\fi

\ifJOURNAL
\documentclass{article}
\usepackage{amsmath,amsfonts,amssymb,amsthm,latexsym,graphicx,algorithm,algorithmic}
\newcommand{\Extra}[1]{}
\fi

\ifWP
\documentclass{article}
\usepackage{amsmath,amsfonts,amssymb,amsthm,latexsym,graphicx,algorithm,algorithmic}

\makeatletter

\newif\iftwodates
\twodatesfalse

\renewcommand\maketitle{\begin{titlepage}%
  \let\footnotesize\small
  \let\footnoterule\relax
  \let \footnote \thanks
  \null\vfil
  \vskip 30\p@
  \begin{center}%
    {\LARGE \bf \@title \par}%
    \vskip 3em%
    {\large
     \lineskip .75em%
     \begin{tabular}[t]{c}%
       \@author
     \end{tabular}\par}%
     \vskip 1.5em%
  \end{center}\par
  \vfill
  \begin{center}
    \raisebox{1.5cm}{\includegraphics[width=0.58\textwidth]%
      {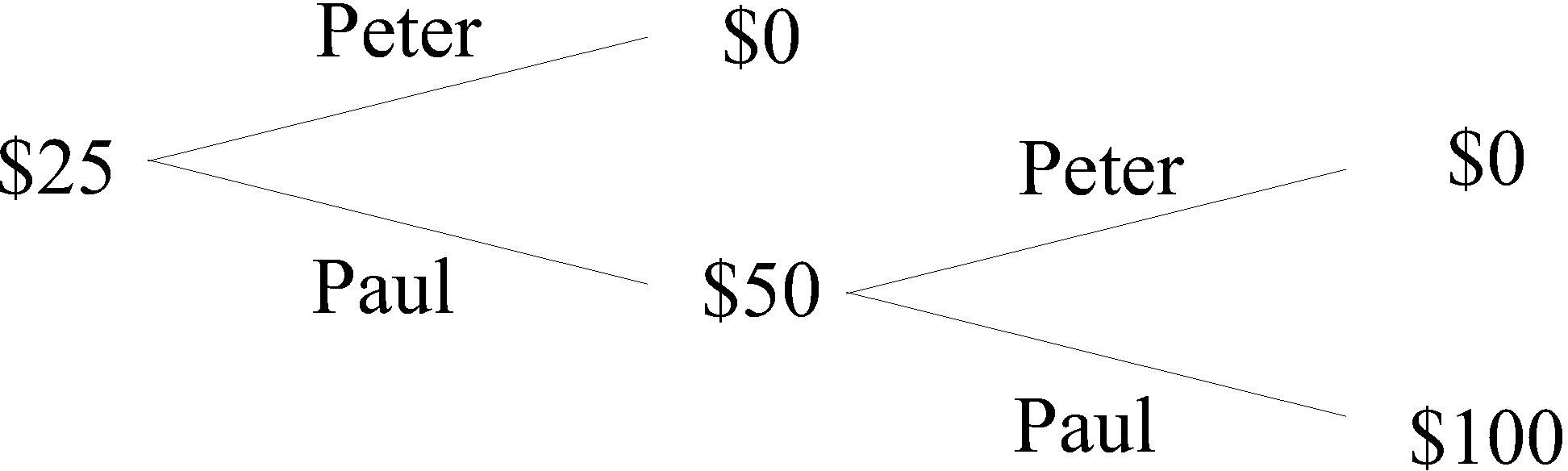}}%
    \hskip 3em%
    \includegraphics[width=0.29\textwidth]%
      {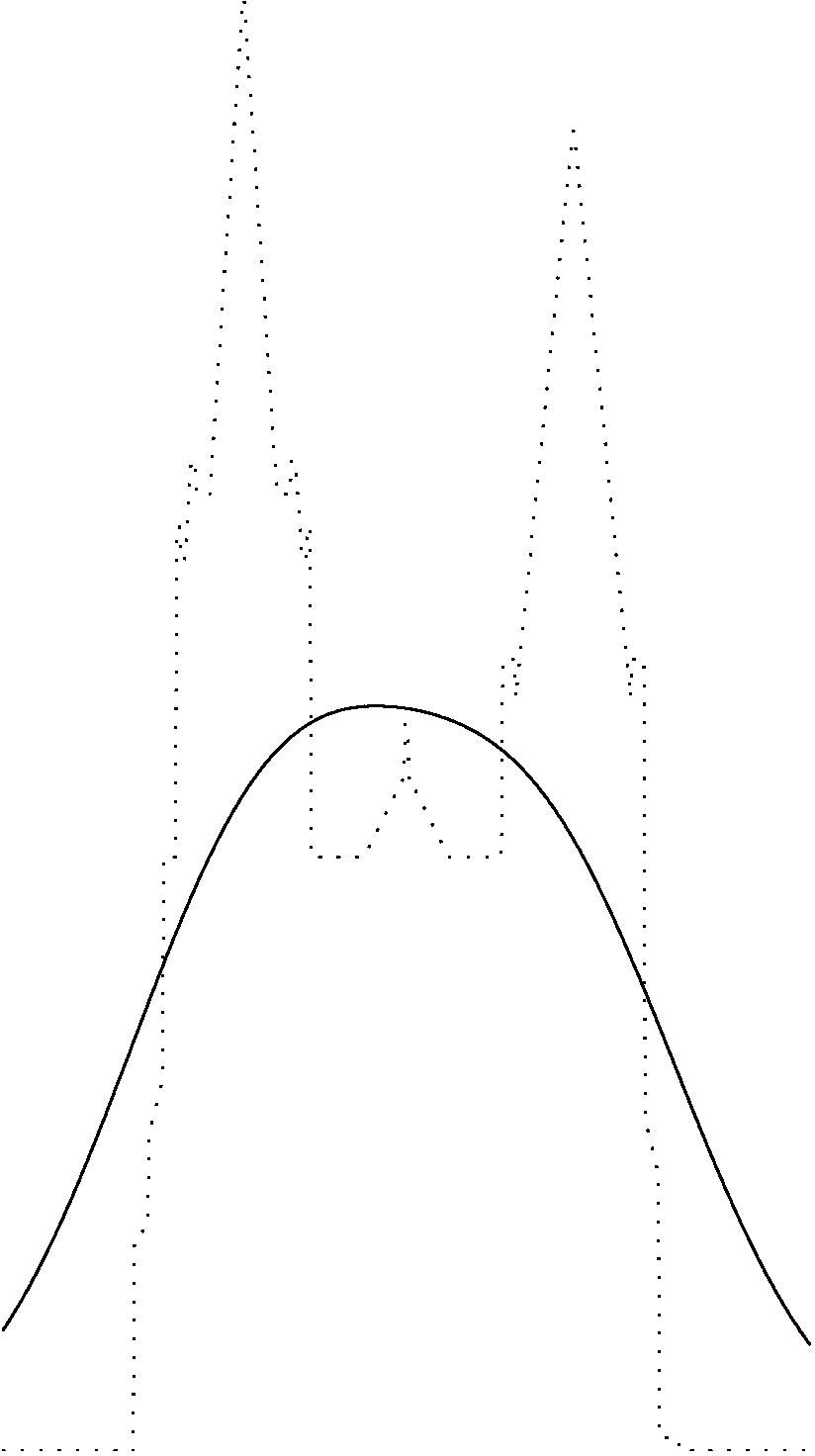}%
  \end{center}
  \@thanks
  \vfill
  \begin{center}
    {\large \bf The Game-Theoretic Probability and Finance Project}
  \end{center}
  \begin{center}
    {\large Working Paper \#\No}
  \end{center}
  \begin{center}
    {\iftwodates\large First posted \firstposted.
    Last revised \@date.\else\large\@date\fi}
  \end{center}
  \begin{center}
    Project web site:\\
    http://www.probabilityandfinance.com
  \end{center}
  \end{titlepage}%
  \setcounter{footnote}{0}%
  \global\let\thanks\relax
  \global\let\maketitle\relax
  \global\let\@thanks\@empty
  \global\let\@author\@empty
  \global\let\@date\@empty
  \global\let\@title\@empty
  \global\let\title\relax
  \global\let\author\relax
  \global\let\date\relax
  \global\let\and\relax
}

\renewenvironment{abstract}{%
  \titlepage
  \null\vfil
  \@beginparpenalty\@lowpenalty
  \begin{center}%
    \Large \bfseries \abstractname
    \@endparpenalty\@M
  \end{center}}%
  {\par\vfill\tableofcontents\endtitlepage}

\renewenvironment{thebibliography}[1]
  {\section*{\refname}%
  \addcontentsline{toc}{section}{\refname}
  \@mkboth{\MakeUppercase\refname}{\MakeUppercase\refname}%
  \list{\@biblabel{\@arabic\c@enumiv}}%
    {\settowidth\labelwidth{\@biblabel{#1}}%
    \leftmargin\labelwidth
    \advance\leftmargin\labelsep
    \@openbib@code
    \usecounter{enumiv}%
    \let\p@enumiv\@empty
    \renewcommand\theenumiv{\@arabic\c@enumiv}}%
    \sloppy
    \clubpenalty4000
    \@clubpenalty \clubpenalty
    \widowpenalty4000%
    \sfcode`\.\@m}
    {\def\@noitemerr
    {\@latex@warning{Empty `thebibliography' environment}}%
  \endlist}

\makeatother

\newcommand{\Extra}[1]{}
\fi

\ifBASIC
\documentclass{article}
\usepackage{amsmath,amsfonts,amssymb,amsthm,latexsym,graphicx,algorithm,algorithmic}
\newcommand{\Extra}[1]{}
\fi

\ifFULL
  \usepackage{color}
  \renewcommand{\Extra}[1]{\blue{#1}}
  
  \newcommand{\blue}[1]{\textcolor{blue}{#1}}
  \newcommand{\bluebegin}{\begingroup\color{blue}}
  \newcommand{\blueend}{\endgroup}

\fi

\allowdisplaybreaks[4]

\newcommand{\Vladimir}{Vladimir}
\newcommand{\DOT}{.}

\ifnotLATIN
  \input{OT2enc.def}

\fi

\gdef\mathenddots{\mathinner{\ldotp\ldotp\ldotp\ldotp}}

\newcommand{\dd}{\mathrm{d}}		

\newcommand{\K}{\mathcal{K}}		
\newcommand{\EEE}{\mathcal{E}}		

\DeclareMathOperator{\III}{\mathbb{I}}		

\newcommand{\bbbn}{\mathbb{N}}                  
\newcommand{\bbbr}{\mathbb{R}}			

\newcommand{\XXX}{\mathcal{X}}		

\newcommand{\bbbrbar}{\overline{\mathbb{R}}}	

\newcommand{\bbbrbarXXX}{\smash{\bbbrbar}\vphantom{\bbbr}^{\XXX}}

\newtheorem{theorem}{Theorem}
\newtheorem*{theorem*}{Theorem}

\newtheorem{lemma}{Lemma}
\newtheorem{corollary}{Corollary}

\theoremstyle{definition}

\makeatletter
\newenvironment{protocol}{%
\renewcommand{\ALG@name}{Protocol}
\begin{algorithm}
}{%
\end{algorithm}
\renewcommand{\ALG@name}{Algorithm}
}
\makeatother

\title{Insuring against loss of evidence in game-theoretic probability}

\ifJOURNAL
\author{A. Philip Dawid$^1$, Steven de Rooij$^1$, Glenn Shafer$^{2,3}$,\\
  Alexander Shen$^4$, Nikolai Vereshchagin$^5$, and Vladimir Vovk$^{3,6}$}
\fi

\ifWP
\author{A. Philip Dawid, Steven de Rooij, Glenn Shafer,\\
  Alexander Shen, Nikolai Vereshchagin, and Vladimir Vovk}
  \newcommand{\No}{34}
\fi

\ifBASIC
\author{A. Philip Dawid, Steven de Rooij, Glenn Shafer,\\
  Alexander Shen, Nikolai Vereshchagin, and Vladimir Vovk}
\fi

\begin{document}
\maketitle

\begin{abstract}
  We consider the game-theoretic scenario of testing the performance
  of Forecaster by Sceptic who gambles against the forecasts.
  Sceptic's current capital is interpreted as the amount of evidence
  he has found against Forecaster.
  Reporting the maximum of Sceptic's capital so far
  exaggerates the evidence.
  We characterize the set of all increasing functions that remove the exaggeration.
  This result can be used for insuring against loss of evidence.
\end{abstract}

\section{Introduction}
\label{sec:introduction}

In game-theoretic probability
(see, e.g., \cite{shafer/vovk:2001})
Sceptic is trying to prove Forecaster wrong
by gambling against him:
the values of Sceptic's capital $\K_n$ measure the changing evidence against Forecaster.
It is always assumed that Sceptic's initial capital is $\K_0=1$,
and Sceptic is required to ensure that $\K_n\ge0$ at each time $n$.

The evidence, however, can be both gained and lost.
When $\K_n$ becomes large at some time $n$,
Forecaster's performance begins to look poor,
but then $\K_i$ for some later time $i$ may be lower
and make Forecaster look better.
Our result will show that Sceptic can avoid losing too much evidence,
with a modest price to pay for this.

Suppose we exaggerate the evidence against Forecaster
by considering not the current value $\K_n$ of his capital
but the greatest value so far: 
$$
  \K^*_n := \max_{i\le n} \K_i.
$$
Continuing research started in \cite{GTP33},
we show that there are many functions $F:[1,\infty)\to[0,\infty)$
such that:
(1) $F(y)\to\infty$ as $y\to\infty$ almost as fast as $y$;
(2) Sceptic's moves can be modified on-line
in such a way that the modified moves lead to capital
\begin{equation}\label{eq:goal}
  \K'_n
  \ge
  F(\K_n^*),
  \quad
  n=1,2,\mathenddots
\end{equation}
Sceptic who is worried about losing evidence can use a middle approach
securing him capital $c\K_n+(1-c)\K'_n$ at each time $n$
for a constant $c\in(0,1)$.
This way he may sacrifice a fraction $1-c$ of his capital
but gets insurance against losing the bulk of his evidence.
See Section~\ref{sec:loss} for details.

Technically,
we characterize the set of increasing functions $F$
for which (\ref{eq:goal}) can be achieved.
In \cite{GTP33}
a similar result is proved
in the framework of measure-theoretic probability.
The latter corresponds to the case where Sceptic's strategy is known in advance
(it involves some other simplifying assumptions,
such as additivity and even $\sigma$-additivity,
but they are less important in our current context).
The situation when Sceptic's strategy is known is much easier,
and \cite{GTP33} uses a simple method based on L\'evy's zero-one law
(see \cite{GTP29} for the game-theoretic version of L\'evy's law).
The method of this article is completely different
and is based on the idea of stopping and combining capital processes.
This idea is known and has been used in, e.g.,
\cite{el-yaniv/etal:2001} (Theorem 1,
based on Leonid Levin's personal communication)
and \cite{shafer/vovk:2001} (Lemma 3.1);
we show that it gives optimal results in our current framework.

In this article the words such as ``positive'' and ``increasing''
will be used in the wide sense of the inequalities $\le$ and $\ge$.
The set of real numbers is $\bbbr$
and the set of natural numbers is $\bbbn:=\{1,2,\ldots\}$.
The extended real line $[-\infty,\infty]$ is denoted $\bbbrbar$,
and we use the convention $\infty+(-\infty):=\infty$.
If $E$ is some property,
$\III_{\{E\}}$ is defined to be $1$ if $E$ is satisfied
and $0$ if not.

\section{Calibrating exaggerated evidence}
\label{sec:calibrating}

Our prediction protocol involves four players:
Forecaster, Sceptic, Rival Sceptic, and Reality.

\begin{protocol}[H]
  \caption{Competitive scepticism}
  \label{prot:competitive}
  \begin{algorithmic}
    \STATE $\K_0:=1$ and $\K'_0:=1$
    \FOR{$n=1,2,\dots$}
      \STATE Forecaster announces $\EEE_n\in\mathbf{E}$
      \STATE Sceptic announces $f_n\in[0,\infty]^{\XXX}$ such that $\EEE_n(f_n)\le\K_{n-1}$
      \STATE Rival Sceptic announces $f'_n\in[0,\infty]^{\XXX}$
        such that $\EEE_n(f'_n)\le\K'_{n-1}$
      \STATE Reality announces $x_n\in\XXX$
      \STATE $\K_n:=f_n(x_n)$ and $\K'_n:=f'_n(x_n)$
    \ENDFOR
  \end{algorithmic}
\end{protocol}

\noindent
The parameter of the protocol is a set $\XXX$,
from which Reality chooses her moves;
$\mathbf{E}$ is the set of all ``outer probability contents'' on $\XXX$
(to be defined momentarily).
We always assume that $\XXX$ contains at least two distinct elements.
The reader who is not interested in the most general statement of our result
can interpret $\mathbf{E}$ as the set of all expectation functionals
$\EEE:f\mapsto\int f \dd P$,
$P$ being a probability measure on a fixed $\sigma$-algebra on $\XXX$;
in this case Sceptic and Rival Sceptic are required to output
functions that are measurable w.r.\ to that $\sigma$-algebra.

In general,
an \emph{outer probability content} on $\XXX$
is a function $\EEE:\bbbrbarXXX\to\bbbrbar$
(where $\bbbrbarXXX$ is the set of all functions $f:\XXX\to\bbbrbar$)
that satisfies the following four axioms:
\begin{enumerate}
\item\label{ax:order}
  If $f,g\in\bbbrbarXXX$ and $f\le g$,
  then $\EEE(f)\le\EEE(g)$.
\item\label{ax:scaling}
  If $f\in\bbbrbarXXX$ and $c\in(0,\infty)$, then $\EEE(cf)=c\EEE(f)$.
\item\label{ax:sum}
  If $f,g\in\bbbrbarXXX$, then $\EEE(f+g)\le\EEE(f)+\EEE(g)$.
\item\label{ax:norm}
  For each $c\in\bbbr$,
  $\EEE(c)=c$,
  where the $c$ in parentheses
  is the function in $\bbbrbarXXX$ that is identically equal to $c$.
\end{enumerate}
An axiom of $\sigma$-subadditivity on $[0,\infty]^{\XXX}$
is sometimes added to this list,
but we do not need it in this article.
(And it is surprising how rarely it is needed in general:
see, e.g., \cite{GTP29}.)
In our terminology we follow \cite{hoffmann-jorgensen:1987} and \cite{GTP29}.
Upper previsions studied in the theory of imprecise probabilities
(see, e.g., \cite{decooman/hermans:2008})
are closely related to (but somewhat more restrictive than)
outer probability contents.

Protocol~\ref{prot:competitive} describes a perfect-information game
in which Sceptic tries to discredit the outer probability contents $\EEE_n$
issued by Forecaster
as a faithful description of how Reality produces $x_n\in\XXX$.
The players make their moves sequentially in the indicated order.
At each step Sceptic and Rival Sceptic choose gambles $f_n$ and $f'_n$
on how $x_n$ is going to come out,
and their resulting capitals are $\K_n$ and $\K'_n$, respectively.
Discarding capital is allowed,
but Sceptic and Rival Sceptic are required to ensure
that $\K_n\ge0$ and $\K'_n\ge0$, respectively;
this is achieved by requiring that $f_n$ and $f'_n$ should be positive.

Let us call an increasing function
$F:[1,\infty)\to[0,\infty)$ a \emph{capital calibrator}
if there exists a strategy for Rival Sceptic that guarantees
$\K'_n\ge F(\K_n^*)$ for all $n$,
with $F(\infty)$ understood to be $\lim_{y\to\infty}F(y)$.
We say that a capital calibrator $F$ \emph{dominates} a capital calibrator $G$
if $F(y)\ge G(y)$ for all $y\in[1,\infty)$.
We say that $F$ \emph{strictly dominates} $G$
if $F$ dominates $G$ and $F(y)>G(y)$ for some $y\in[1,\infty)$.
A capital calibrator is \emph{admissible} if it is not strictly dominated
by any other capital calibrator.
\begin{theorem}\label{thm:calibration}
\begin{enumerate}
\item
  An increasing function $F:[1,\infty)\to[0,\infty)$ is a capital calibrator
  if and only if
  \begin{equation}\label{eq:reducer-le}
    \int_1^{\infty}
    \frac{F(y)}{y^2}
    \dd y
    \le
    1.
  \end{equation}
\item
  Any capital calibrator is dominated by an admissible capital calibrator.
\item
  A capital calibrator is admissible if and only if
  it is right-continuous and
  \begin{equation}\label{eq:reducer-eq}
    \int_1^{\infty}
    \frac{F(y)}{y^2}
    \dd y
    =
    1.
  \end{equation}
\end{enumerate}
\end{theorem}
\begin{proof}
  First we prove that any increasing function $F:[1,\infty)\to[0,\infty)$ satisfying
  \begin{equation}\label{eq:alternative}
    F(y)
    =
    \int_{[1,y]}
    u
    P(\dd u),
    \quad
    \forall y\in[1,\infty),
  \end{equation}
  for a probability measure $P$ on $[1,\infty)$
  is a capital calibrator.
  For each $u\ge1$, define the following strategy for Rival Sceptic:
  at step $n$, the strategy outputs
  $$
    f_n^{(u)}
    :=
    \begin{cases}
      f_n & \text{if $\K^*_{n-1}<u$}\\
      u & \text{otherwise}
    \end{cases}
  $$
  as Rival Sceptic's move $f'_n$.
  Let us check that this is a valid strategy,
  i.e., that $\EEE_n(f_n^{(u)})\le\K^{(u)}_{n-1}$, $n\in\bbbn$,
  where $\K^{(u)}$ is defined by $\K^{(u)}_0:=1$
  and $\K^{(u)}_n:=f_n^{(u)}(x_n)$ for $n\in\bbbn$.
  There are three cases to consider:
  \begin{enumerate}
  \item
    If $\K^*_{n-1}<u$,
    we have $\K^{(u)}_{n-1}=\K_{n-1}$ and
    $
      \EEE_n(f_n^{(u)})=\EEE_n(f_n)\le\K_{n-1}=\K^{(u)}_{n-1}
    $.
  \item
    If $n$ is the smallest number for which $\K^*_{n-1}\ge u$,
    we have $\K^{(u)}_{n-1}=\K_{n-1}\ge u$ and
    $
      \EEE_n(f_n^{(u)})=\EEE_n(u)=u\le\K^{(u)}_{n-1}
    $.
  \item
    Otherwise, we have $\K^{(u)}_{n-1}=u$ and so
    $
      \EEE_n(f_n^{(u)})=\EEE_n(u)=u=\K^{(u)}_{n-1}
    $.
  \end{enumerate}
  Set $f'_n(x):=\int_{[1,\infty)} f_n^{(u)}(x) P(\dd u)$,
  $x\in\XXX$;
  this gives $\K'_n=\int_{[1,\infty)}\K_n^{(u)} P(\dd u)$
  when we set $x$ to $x_n$.
  Let us check that this is a valid strategy for Rival Sceptic,
  i.e., that $\EEE_n(f'_n)\le\K'_{n-1}$ for all $n\in\bbbn$.
  This is now obvious if $\EEE_n$ are expectation functionals,
  and in general we have
  \begin{align*}
    \EEE_n(f'_n)
    &=
    \EEE_n
    \left(
      \int_{[1,\infty)}
      f_n^{(u)}
      P(\dd u)
    \right)\\
    &=
    \EEE_n
    \left(
      \int_{[1,\infty)}
      \left(
        \III_{\{\K^*_{n-1}<u\}} f_n
        +
        \III_{\{\K^*_{n-1}\ge u\}} u
      \right)
      P(\dd u)
    \right)\\
    &=
    \EEE_n
    \left(
      P((\K^*_{n-1},\infty)) f_n
      +
      \int_{[1,\K^*_{n-1}]} u P(\dd u)
    \right)\\
    &\le
    P((\K^*_{n-1},\infty)) \K_{n-1}
    +
    \int_{[1,\K^*_{n-1}]} u P(\dd u)\\
    &=
    \int_{(\K^*_{n-1},\infty)} \K_{n-1} P(\dd u)
    +
    \int_{(\K^*_{n-2},\K^*_{n-1}]} u P(\dd u)
    +
    \int_{[1,\K^*_{n-2}]} u P(\dd u)\\
    &\le
    \int_{(\K^*_{n-1},\infty)} \K_{n-1}^{(u)} P(\dd u)
    +
    \int_{(\K^*_{n-2},\K^*_{n-1}]} \K_{n-1}^{(u)} P(\dd u)
    +
    \int_{[1,\K^*_{n-2}]} \K_{n-1}^{(u)} P(\dd u)\\
    &=
    \int_{[1,\infty)}
    \K_{n-1}^{(u)}
    P(\dd u)
    =
    \K'_{n-1}.
  \end{align*}
  The last inequality used the analysis of the three cases above.
  For small values of $n$, our convention was $\K^*_0:=1$ and $\K^*_{-1}:=1$.
  Notice that our argument only used Axioms~\ref{ax:scaling}--\ref{ax:norm}
  for outer probability contents;
  no $\sigma$-subadditivity was required.
  This strategy will guarantee
  \begin{equation}\label{eq:guarantee}
    \K'_n
    =
    \int_{[1,\infty)} \K_n^{(u)} P(\dd u)
    \ge
    \int_{[1,\K_n^*]} \K_n^{(u)} P(\dd u)
    \ge
    \int_{[1,\K_n^*]} u P(\dd u)
    =
    F(\K_n^*).
  \end{equation}

  We can now finish the proof of the statement ``if'' in part~1 of the theorem,
  which says that any increasing function $F:[1,\infty)\to[0,\infty)$
  satisfying~(\ref{eq:reducer-le}) is a capital calibrator.
  Without loss of generality we can 
  assume that $F$ is right-continuous and that (\ref{eq:reducer-eq}) holds.
  It remains to apply Lemma~\ref{lem:reducer} below.

  Let us now check that every capital calibrator satisfies~(\ref{eq:reducer-le}).
  Suppose a capital calibrator $F$ violates~(\ref{eq:reducer-le}).
  We can decrease $F$ so that,
  for some $a>1$ and $N\in\bbbn$,
  it is constant in each interval $[a^{n-1},a^n)$, $n=1,\ldots,N$,
  is zero in $[a^N,\infty)$,
  and still violates~(\ref{eq:reducer-le}).
  Of course, $F$ is still a capital calibrator.
  The substitution $x=1/y$ shows that
  $
    \int_0^1 F(1/x) \dd x
    >
    1
  $,
  which can be rewritten as
  \begin{equation}\label{eq:sum}
    F(1)
    \left(
      1 - \frac1a
    \right)
    +
    F(a)
    \left(
      \frac1a - \frac{1}{a^2}
    \right)
    +\cdots+
    F(a^{N-1})
    \left(
      \frac{1}{a^{N-1}} - \frac{1}{a^N}
    \right)
    >
    1.
  \end{equation}
  Suppose, without loss of generality, that $\XXX\supseteq\{0,1\}$,
  and let Forecaster always choose
  $$
    \EEE_n(f)
    :=
    \frac1a
    f(1)
    +
    \left(
      1 - \frac1a
    \right)
    f(0),
    \quad
    n\in\bbbn.
  $$
  Let Sceptic play the strategy of always betting all his capital on $1$:
  $f_n(1):=a\K_{n-1}$ and $f_n(x):=0$ for $x\ne1$.
  Then $\K^*_N=a^n$ where $n$ is the number of $1$s output by Reality
  before the first element different from $1$
  (except that $n=N$ if Reality outputs only $1$s
  during the first $N$ steps).
  Backward induction shows that the initial capital $\K'_0$
  required to ensure $\K'_N\ge F(\K^*_N)$
  must be at least
  \begin{multline*}
    F(a^N)
    \left(
      \frac1a
    \right)^N
    +
    F(a^{N-1})
    \left(
      \frac1a
    \right)^{N-1}
    \left(
      1 - \frac1a
    \right)
    +
    F(a^{N-2})
    \left(
      \frac1a
    \right)^{N-2}
    \left(
      1 - \frac1a
    \right)
    \\+\cdots+
    F(a)
    \frac1a
    \left(
      1 - \frac1a
    \right)
    +
    F(1)
    \left(
      1 - \frac1a
    \right)
    >
    1;
  \end{multline*}
  the inequality follows from~(\ref{eq:sum}),
  but we know that it is false as $\K'_0=1$.

  We have proved part~1 of the theorem.
  Part~3 is now obvious, and part~2 follows from parts~1 and~3.
\end{proof}

The following lemma was used in the proof of Theorem~\ref{thm:calibration}.
\begin{lemma}\label{lem:reducer}
  An increasing right-continuous function $F:[1,\infty)\to [0,\infty)$
  satisfies~(\ref{eq:reducer-eq}) if and only if (\ref{eq:alternative})
  holds for some probability measure $P$ on $[1,\infty)$.
\end{lemma}
\begin{proof}
  Let us first check that the existence of a probability measure $P$
  satisfying~(\ref{eq:alternative})
  implies~(\ref{eq:reducer-eq}).
  We have:
  \begin{multline}
    \int_{[1,\infty)}
    \frac{F(y)}{y^2}
    \dd y
    =
    \int_{[1,\infty)}
    \int_{[1,y]}
    \frac{u}{y^2}
    P(\dd u)
    \dd y\\
    =
    \int_{[1,\infty)}
    \int_{[u,\infty)}
    \frac{u}{y^2}
    \dd y
    P(\dd u)
    =
    \int_{[1,\infty)}
    P(\dd u)
    =
    1.
    \label{eq:integral}
  \end{multline}

  It remains to check that any increasing right-continuous
  $F:[1,\infty)\to[0,\infty)$
  satisfying~(\ref{eq:reducer-eq}) satisfies~(\ref{eq:alternative})
  for some probability measure $P$ on $[1,\infty)$.
  Let $Q$ be the measure on $[1,\infty)$
  ($\sigma$-finite but not necessarily a probability measure)
  with distribution function $F$,
  in the sense that $Q([1,y])=F(y)$
  for all $y\in[1,\infty)$.
  Set $P(\dd u):=(1/u)Q(\dd u)$.
  We then have~(\ref{eq:alternative}),
  and the calculation (\ref{eq:integral})
  shows that the $\sigma$-finite measure $P$ must be a probability measure
  (were it not, we would not have an equality in~(\ref{eq:reducer-eq})).
\end{proof}

According to~(\ref{eq:reducer-eq}),
the functions
\begin{equation}\label{eq:class}
  F(y)
  :=
  \alpha y^{1-\alpha}
\end{equation}
are admissible capital calibrators for any $\alpha\in(0,1)$.

\section{Insuring against loss of evidence}
\label{sec:loss}

As we saw in Section~\ref{sec:introduction},
there is a simple way to use Theorem~\ref{thm:calibration}
for insuring against loss of evidence.
The following corollary says that it leads to an optimal result.

\begin{corollary}\label{cor:insurance}
  Let $c\ge0$ and $F:[1,\infty)\to[0,\infty)$ be an increasing function.
  Rival Sceptic has a strategy ensuring
  \begin{equation}\label{eq:insurance}
    \K'_n\ge c\K_n+F(\K_n^*)
  \end{equation}
  if and only if $c$ and $F$ satisfy
  \begin{equation}\label{eq:condition}
    \int_1^{\infty} \frac{F(y)}{y^2} \dd y \le 1-c.
  \end{equation}
\end{corollary}
\begin{proof}
  Suppose (\ref{eq:condition}) is satisfied;
  in particular, $c\in[0,1]$.
  Using $cf_n+(1-c)f'_n$ as Rival Sceptic's strategy,
  where $f_n$ are Sceptic's moves and $f'_n$ are Rival Sceptic's moves
  guaranteeing $\K'_n\ge\frac{1}{1-c}F(\K^*_n)$ (cf.\ Theorem~\ref{thm:calibration}),
  we can see that Rival Sceptic can guarantee (\ref{eq:insurance}).

  Now suppose Rival Sceptic can ensure (\ref{eq:insurance}),
  but (\ref{eq:condition}) is violated.
  As in the proof of Theorem~\ref{thm:calibration},
  we can decrease $F$ so that,
  for some $a>1$ and $N\in\bbbn$,
  it is constant in each interval $[a^{n-1},a^n)$, $n=1,\ldots,N$,
  is zero in $[a^N,\infty)$,
  and still violates~(\ref{eq:condition}).
  Similarly to (\ref{eq:sum}), we have
  \begin{equation*}
    F(1)
    \left(
      1 - \frac1a
    \right)
    +
    F(a)
    \left(
      \frac1a - \frac{1}{a^2}
    \right)
    +\cdots+
    F(a^{N-1})
    \left(
      \frac{1}{a^{N-1}} - \frac{1}{a^N}
    \right)
    >
    1-c.
  \end{equation*}
  Suppose $\XXX\supseteq\{0,1\}$
  and define Forecaster's and Sceptic's strategies as before.
  Now backward induction shows that the initial capital $\K'_0$
  required to ensure $\K'_N\ge c\K_N+F(\K^*_N)$
  must be at least
  \begin{multline*}
    c a^N
    \left(
      \frac1a
    \right)^N
    +
    F(a^N)
    \left(
      \frac1a
    \right)^N
    +
    F(a^{N-1})
    \left(
      \frac1a
    \right)^{N-1}
    \left(
      1 - \frac1a
    \right)\\
    +
    F(a^{N-2})
    \left(
      \frac1a
    \right)^{N-2}
    \left(
      1 - \frac1a
    \right)
    +\cdots+
    F(a)
    \frac1a
    \left(
      1 - \frac1a
    \right)
    +
    F(1)
    \left(
      1 - \frac1a
    \right)\\
    >
    c + (1-c) = 1.
  \end{multline*}
  This contradicts $\K'_0=1$.
\end{proof}

According to (\ref{eq:class}) and (\ref{eq:insurance}),
Rival Sceptic can guarantee
\begin{equation}\label{eq:alpha-insurance}
  \K'_n
  \ge
  c\K_n
  +
  (1-c)\alpha
  (\K^*_n)^{1-\alpha}
\end{equation}
for any constants $c\in[0,1]$ and $\alpha\in(0,1)$.

Corollary~\ref{cor:insurance} does not mean that (\ref{eq:alpha-insurance})
or, more generally, (\ref{eq:insurance}) cannot be improved;
it only says that the improvement will not be significant enough
to decrease the coefficient in front of $\K_n$.
For example, if we do not discard the term
$\int_{(\K_n^*,\infty)}\K_n^{(u)}P(\dd u)$
in (\ref{eq:guarantee}),
we will obtain
\begin{equation}\label{eq:stronger-guarantee}
  \K'_n
  \ge
  P((\K_n^*,\infty))
  \K_n
  +
  F(\K_n^*).
\end{equation}
The coefficient $P((\K_n^*,\infty))$ in front of $\K_n$ tends to zero
as $\K_n^*\to\infty$.

In particular,
using (\ref{eq:stronger-guarantee}) allows us to improve~(\ref{eq:alpha-insurance})
to
\begin{equation*}
  \K'_n
  \ge
  c\K_n
  +
  (1-c)(1-\alpha)(\K_n^*)^{-\alpha}\K_n
  +
  (1-c)\alpha(\K^*_n)^{1-\alpha}.
\end{equation*}

\end{document}